\documentclass[12pt,oneside]{article}

\usepackage{latexsym,amssymb,amsfonts,amsmath,amsthm}

\usepackage{mathrsfs}

\allowdisplaybreaks[3]

\theoremstyle{plain}
\newtheorem{thm}{Theorem}[section]

\newtheorem{lem}[thm]{Lemma}
\newtheorem{prop}[thm]{Proposition}

\theoremstyle{definition}
\newtheorem*{notation}{Notation}

\theoremstyle{remark}
\newtheorem*{rems}{Remark}

\newcommand{\vtsp}{\hspace{0.1em}}

\newlength{\phantomheight}

\DeclareMathOperator{\Trace}{Tr}

\newcommand{\R}{\mathbb{R}}
\newcommand{\N}{\mathbb{N}}
\newcommand{\NZ}{\mathbb{N}_0}
\newcommand{\C}{\mathbb{C}}
\newcommand{\Z}{\mathbb{Z}}
\renewcommand{\epsilon}{\varepsilon}

\newcommand{\ipd}[2]{\langle{#1},{#2}\rangle}
\newcommand{\bigipd}[2]{\bigl\langle{#1},\vtsp{#2}\bigr\rangle}

\newcommand{\abs}[1]{\lvert{#1}\rvert}
\newcommand{\lrabs}[1]{\left\lvert{#1}\right\rvert}
\newcommand{\bigabs}[1]{\bigl\lvert{#1}\bigr\rvert}
\newcommand{\norm}[1]{\lVert{#1}\rVert}

\renewcommand{\d}{\,d}

\DeclareMathOperator{\Cos}{Cos}

\newcommand{\op}{H}
\newcommand{\pot}{V}
\newcommand{\res}{R}
\newcommand{\zqlat}{\Lambda}
\newcommand{\qlat}{\Lambda'}
\newcommand{\md}{\gamma}
\renewcommand{\a}{\mathbf{a}}
\newcommand{\potcpt}[1][\a]{\pot_{#1}}
\newcommand{\Uop}[1][\a]{U_{#1}}
\newcommand{\Ucpt}[2][\a]{U_{#1}^{#2}}
\newcommand{\TR}{{T^*\R}}
\newcommand{\efn}[1]{\phi_{#1}}

\newcommand{\aba}[1]{\abs{#1}_{\alpha}}
\newcommand{\He}{\mathcal H}

\title{Asymptotics for the Eigenvalues of the Harmonic Oscillator with
  a Quasi-Periodic Perturbation}

\author{Daniel M.~Elton}

\begin{document}

\maketitle

\begin{abstract}
We consider operators of the form $\op+\pot$ where $\op$ is the
one-dimensional harmonic oscillator and $\pot$ is a zero-order
pseudo-differential operator which is quasi-periodic in an appropriate
sense (one can take $\pot$ to be multiplication by a periodic function
for example). It is shown that the eigenvalues of $\op+\pot$ have
asymptotics of the form $\lambda_n(\op+\pot)=\lambda_n(\op)+W(\sqrt
n)n^{-1/4}+O(n^{-1/2}\ln(n))$ as $n\to+\infty$, 
where $W$ is a quasi-periodic function
which can be defined explicitly in terms of $\pot$. 
\end{abstract}

\section{Introduction}

The one-dimensional harmonic oscillator is the operator
\[
\op\,=\,-\frac{d^2}{dx^2}+(\alpha x)^2,
\]
where $\alpha$ is a positive parameter. We can consider $\op$ as an
unbounded self-adjoint operator acting on $L^2(\R)$. The determination
of the spectrum of $\op$ is a classical problem --- 
virtually any introductory book on quantum mechanics
has a section devoted to this topic. In particular $\op$
has a compact resolvent and hence a discrete spectrum. Furthermore, the
eigenvalues of $\op$ are simple and can be enumerated as  
\[
\lambda_n(\op)=\alpha(2n+1),\qquad n\in\NZ.
\]
A normalised eigenfunction corresponding to $\lambda_n(\op)$ can be
chosen as
\begin{equation}
\label{defofefn:eq}
\efn{n}(x)\,=\,\frac{\alpha^{1/4}}{\sqrt{n!\vtsp 2^n\sqrt\pi}}
\,e^{-\alpha x^2/2}\,\He_n(\sqrt\alpha x),
\end{equation}
where $\He_n$ is the $n$--th Hermite polynomial.

\medskip

The purpose of this paper is to study the large $n$ asymptotics of the
eigenvalues of the perturbed operator $\op+\pot$ when $\pot$ is a
self-adjoint quasi-periodic pseudo-differential operator of order
$0$. More precisely, we assume $V$ can be written in the form
\begin{equation}
\label{defofpot:eq}
\pot\,=\,\sum_{\a\in\zqlat}\,\potcpt\Uop
\end{equation}
where $\zqlat\subset\TR\cong\R^2$ is a countable discrete index set
and, for each $\a=(a_x,a_\xi)\in\TR$, we define $\Uop$ to be the
unitary operator on $L^2(\R)$ given by
\begin{equation}
\label{defofUop:eq}
\Uop\phi(x)\,=\,e^{ia_xa_\xi/2}e^{ia_xx}\phi(x+a_\xi).
\end{equation}
The $\potcpt$'s are just complex coefficients. 

Since $\Uop^*=\Uop[-\a]$ for any $\a\in\TR$, the condition that $\pot$
is self-adjoint can be rewritten as the requirement
\[
\a\in\zqlat\;\Longrightarrow\;-\a\in\zqlat
\quad\ \text{and}\ \quad
\potcpt[-\a]\,=\,\overline{\potcpt},\quad\a\in\zqlat.
\]
We will also assume the $\potcpt$'s satisfy the following condition
(essentially a regularity assumption);
\begin{equation}
\label{potcond:eq}
\sum_{\a\in\zqlat}\,\abs{\a}^3\abs{\potcpt}
\,<\,+\infty.
\end{equation}
In particular, this condition ensures that the right hand side of
\eqref{defofpot:eq} is absolutely convergent in operator norm, making
$\pot$ a well defined bounded operator. Since $\op$ has a compact
resolvent the same must then be true for $\op+\pot$; it follows that the
spectrum of $\op+\pot$ also consists of discrete eigenvalues.

\begin{rems}
If we take $\zqlat=\{(\omega m,0)\,|\,m\in\Z\}$ then $\pot$ is the
operator of multiplication by a function with period $\omega$ whose
$m$-th Fourier coefficient is simply $\omega^{1/2}\potcpt[(\omega
m,0)]$. Condition \eqref{potcond:eq} becomes a standard regularity
requirement (that the function $V$ should be a ``bit more'' than $C^3$).

In general we may consider $V$ to be a zero-order
pseudo-differential operator 
with Weyl-symbol $\sum_{\a\in\zqlat}\potcpt e^{i(a_xx+a_\xi\xi)}$
(n.b., $\Uop$ is the operator with Weyl-symbol
$e^{i(a_xx+a_\xi\xi)}$). If $\zqlat$ is a rational periodic lattice
then $V$ will be a periodic operator (in the sense that it commutes
with a specific translation operator). Taking $\zqlat$ to be an
irrational periodic lattice, or an irregular discrete set, leads to a
generalisation of such periodic operators; when we apply
``quasi-periodic'' to $\pot$ we mean this particular type of
generalisation. 
\end{rems}

\medskip

If $\mathbf 0\in\zqlat$ then the corresponding term in $\pot$ is
$\potcpt[\mathbf 0]$ times the identity operator and will thus cause a
simple shift in the spectrum of $\op$ by $\potcpt[\mathbf 0]$. This
term is included in the statement of the main result (Theorem
\ref{mainres1:thm} below) but thereafter we shall assume
$\potcpt[\mathbf 0]=0$. We also set 
$\qlat=\zqlat\!\setminus\!\{\mathbf 0\}$; since $\zqlat$ is discrete,
$\TR\!\setminus\!\qlat$ contains a neighbourhood of $\mathbf 0$.

Define a metric $\aba\cdot$ on $\TR$ by
$\aba\a=(\alpha^{-1}a_x^2+\alpha\vtsp a_\xi^2)^{1/2}$.
This metric is equivalent to the usual metric $\abs\cdot$ so condition
\eqref{potcond:eq} can be rewritten as
\begin{equation}
\label{potcond1:eq}
\sum_{\a\in\qlat}\,\aba\a^p\abs{\potcpt}\,<\,+\infty
\quad\text{for all $p\le3$.}
\end{equation}
The main result of the paper is the following.

\begin{thm}
\label{mainres1:thm}
Suppose $\pot$ given by \eqref{defofpot:eq} satisfies
\eqref{potcond:eq} (or equivalently \eqref{potcond1:eq}).
Then the eigenvalues of the operator $\op+\pot$ satisfy
\[
\lambda_n(\op+\pot)\;=\;\alpha(2n\!+\!1)\,+\,\potcpt[\mathbf 0]\,+\,
W(\sqrt{n})\vtsp n^{-1/4}\;+\;O(n^{-1/2}\ln(n))
\]
as $n\to\infty$, where $W:\R\to\R$ is the quasi-periodic function
defined by 
\begin{equation}
\label{defofauxfnU:eq}
W(\lambda)\,=\,\frac{2^{1/4}}{\sqrt\pi}
\sum_{\a\in\qlat}\,\potcpt\aba\a^{-1/2}
\cos\Bigl(\sqrt2\vtsp\aba\a\lambda-\frac\pi4\Bigr).
\end{equation}
\end{thm}

The presence of the quasi-periodic function $W$ means the first order
asymptotics given by Theorem \ref{mainres1:thm} contain considerably
more information about the operator $\pot$ than one might expect
(c.f.\ the simple power type asymptotics
for the case when $\pot$ is given as multiplication by an element of
$C^\infty_0$ (\cite{PS}) or for the operator
$-d^2/d\theta^2\!+\!V(\theta)$ on $S^1$ (see Theorem 4.2 in \cite{MO})). 
In particular we 
note that if $\pot$ is given as multiplication by a periodic
function, knowledge of the first order asymptotics of
$\lambda_n(\op+\pot)$ allows the Fourier coefficients of
$V$ to be ``half'' determined (the values of
$\potcpt[(-m\omega,0)]+\potcpt[(m\omega,0)]$, $m\in\N$, can be
determined from $W$). 

It is likely that there exists a full asymptotic expansion for
$\lambda_n(\op+\pot)$, involving further terms with quasi-periodic
functions multiplying increasingly negative powers of $n$. Judging by
numerical evidence (for example with the potential
$V(x)=\cos(x)$) the second term in the asymptotics is $O(n^{-3/4})$.
This order (even as an improvement of the remainder estimate in
Theorem \ref{mainres1:thm}) appears to involve reasonable subtle
cancellation effects within the series giving the second term of the
asymptotics; no attempt to deal with this analysis is made here.

\smallskip

\begin{rems}
With an obvious modification to the definition of $W$ and a remainder
estimate of $O(n^{-1/3}\ln(n))$, Theorem \ref{mainres1:thm}
also holds for operators $\pot$ of the form
\[
\pot\,=\,\int_\TR \potcpt\Uop\d^2\a
\quad\text{where $\potcpt$ satisfies}\quad
\int _\TR(\aba\a^{-3/2}+\aba\a^3)\abs\potcpt\d^2\a
\,<\,+\infty.
\]
In this case $\pot$ is a pseudo-differential operator of order zero
whose Weyl-symbol has Fourier transform $2\pi\potcpt$. The $\aba\a^3$
term in the condition on $\potcpt$ is then a regularity condition,
while the $\aba\a^{-3/2}$ term is a generalisation of
quasi-periodicity.
\end{rems}

\smallskip

The proof of Theorem \ref{mainres1:thm} is given in Section \ref{mainrespf:sec}
using standard ideas to express the eigenvalues of $\op+\pot$ in terms
of a series involving the resolvent of $\op$ and the operator
$\pot$. The non-triviality of Theorem \ref{mainres1:thm} is contained in
technical results used to establish the convergence of these
series. These results are obtained in Sections \ref{matres:sec} and
\ref{resest:sec}; estimates for the elements
$\ipd{\pot\efn{k}}{\efn{k'}}$ of the matrix of $\pot$ with respect to the
eigenbasis $\{\efn{k}\,|\,k\in\NZ\}$ are obtained in the former and are
then combined to give resolvent estimates in the latter.

\begin{notation}
We use $C$ to denote any positive real constant whose exact value is
not important but which may depend only on the things it is allowed to
in a given problem. Appropriate function type notation is used in
places to make this clearer whilst subscripts are added if we need to keep track of
the value of a particular constant (e.g.~$C_1(V)$ etc.). 

We use $\norm{T}$, $\norm{T}_1$ and $\norm{T}_2$ to denote the
operator, trace class and Hilbert-Schmidt norms of the operator $T$
respectively.
\end{notation}

\section{Estimates for Matrix Elements}
\label{matres:sec}

The aim of this section is to obtain the necessary estimates for the
matrix elements $\ipd{\pot\efn{k}}{\efn{k'}}$ for all $k,k'\in\NZ$. In
turn these will be estimated via 
\begin{equation}
\label{defofUcpt:eq}
\Ucpt{k,k'}\,:=\,\ipd{\Uop\efn{k}}{\efn{k'}}
\end{equation}
defined for all $\a\in\TR$ and $k,k'\in\NZ$. 
Since the operator $\Uop$ is unitary we immediately get 
\begin{equation}
\label{UisUso1:eq}
\abs{\Ucpt{k,k'}}\,\le\,1.
\end{equation}
To obtain more precise estimates we can use the following 
special function identity (see 7.377 on page 844 of \cite{GRJ})
to find an explicit formula for $\Ucpt{k,k'}$; 
for any $0\le k\le k'$ and $y,z\in\C$ we have 
\begin{equation}
\label{hermiteidGRJ:eq}
\int_\R e^{-x^2}\vtsp\He_{k}(x+y)\vtsp\He_{k'}(x+z)\d x
\;=\;2^{k'}\sqrt\pi\vtsp k!\vtsp z^{k'-k}L_k^{(k'-k)}(-2yz),
\end{equation}
where $L_k^{(k'-k)}$ is the generalised Laguerre polynomial.

\begin{lem}
\label{stanmatele:lem}
For any $0\le k\le k'$ and $\a\in\TR\!\setminus\!\{\mathbf 0\}$ we have
\[
\Ucpt{k,k'}\,=\,\sqrt{\frac{k!}{k'!}}
(\sqrt2\rho e^{i\theta})^{k'-k}\vtsp e^{-\rho^2}\vtsp L_k^{(k'-k)}(2\rho^2)
\]
for some $\theta\in\R$, where
\begin{equation}
\label{defofrho:eq}
\rho\,=\,\frac12\Bigl(\frac{a_x^2}\alpha +\alpha a_\xi^2\Bigr)^{1/2}
\,=\,\frac12\aba\a.
\end{equation}
\end{lem}

\begin{proof}
Introduce the complex number
\[
\omega\,=\,\frac{\sqrt\alpha\vtsp a_\xi}2
-i\,\frac{a_x}{2\sqrt\alpha}.
\]
From \eqref{defofUcpt:eq}, \eqref{defofUop:eq} and \eqref{defofefn:eq} we get
\begin{eqnarray*}
\Ucpt{k,k'}&=&\ipd{\Uop\efn{k}}{\efn{k'}}\\
&=&\frac{\sqrt\alpha\vtsp2^{-(k+k')/2}}{\sqrt{k!k'!\pi}}
e^{ia_xa_\xi/2}\int_\R e^{ia_xx}\vtsp e^{-\alpha(x+a_\xi)^2/2}
\vtsp e^{-\alpha x^2/2}\\[-6pt]
&&\hspace*{6cm}\He_{k}(\sqrt\alpha(x\!+\!a_\xi))
\,\He_{k'}(\sqrt\alpha\vtsp x)\d x\\[4pt]
&=&\frac{2^{-(k+k')/2}}{\sqrt{k!k'!\pi}}
\vtsp e^{\omega^2-\alpha a_\xi^2/2+ia_xa_\xi/2}
\int_\R e^{-x^2}\,\He_{k}(x-\omega+\sqrt\alpha a_\xi)
\,\He_{k'}(x-\omega)\d x\\
&=&\sqrt{\frac{k!}{k'!}}\vtsp 2^{(k'-k)/2}
(-\omega)^{k'-k}e^{\omega^2-\alpha a_\xi^2/2+ia_xa_\xi/2}
L_k^{(k'-k)}\bigl(-2\omega(\omega\!-\!\sqrt\alpha\vtsp a_\xi)\bigr)
\end{eqnarray*}
where the last line follows from \eqref{hermiteidGRJ:eq}.
Now $\abs\omega=\rho$ while
\[
\omega^2-\frac{\alpha a_\xi^2}2+\frac{ia_xa_\xi}2
\;=\;\frac{\alpha a_\xi^2}4-\frac{a_x^2}{4\alpha}
-\frac{\alpha a_\xi^2}{2}-\frac{ia_xa_\xi}2+\frac{ia_xa_\xi}2
\;=\;-\abs\omega^2
\]
and
\[
-2\omega(\omega\!-\!\sqrt\alpha\vtsp a_\xi)
\,=\,-2\omega(-\overline\omega)
\,=\,2\abs{\omega}^2.
\]
The result follows.
\end{proof}

Throughout the remainder of this section we will assume
$\a\in\TR\!\setminus\!\{\mathbf 0\}$ is
fixed and $\rho>0$ is given by \eqref{defofrho:eq}.

\medskip

Laguerre polynomials can be expressed in terms of the confluent
hypergeometric function; using 22.5.54 in \cite{AS} we get
\[
L_k^{(k'-k)}(2\rho^2)
\,=\,\binom{k'}{k}\vtsp M(-k,k'-k+1,2\rho^2).
\]
The confluent hypergeometric function can, in turn, be written as a
pointwise absolutely convergent series of Bessel functions; from
13.3.7 in \cite{AS} we get 
\begin{eqnarray*}
M(-k,k'-k+1,2\rho^2)
&\!=\!&(k'-k)!\,e^{\rho^2}\bigl(\rho^2(k'+k+1)\bigr)^{-(k'-k)/2}\\
&&\quad\sum_{j=0}^\infty A_j\left(\frac{\rho}{(k'+k+1)^{1/2}}\right)^j
J_{k'-k+j}\bigl(2\rho\sqrt{k'+k+1}\bigr),
\end{eqnarray*}
where 
\begin{equation}
\label{defAja:eq}
A_0=1,\quad A_1=0,\quad A_2=\frac12(k'-k+1)
\end{equation}
and, for $j\ge 2$,
\begin{equation}
\label{defAjb:eq}
(j+1)A_{j+1}
\,=\,(j+k'-k)A_{j-1}-(k'+k+1)A_{j-2}.
\end{equation}
It follows from Lemma \ref{stanmatele:lem} that 
\begin{equation}
\label{efeimxpkpn:eq}
\Ucpt{k,k'}\;=\;e^{i(k'-k)\theta}\sqrt{\smash[b]{\big.F_{k',k}}}\,
\sum_{j=0}^\infty A_j\left(\frac{\rho}{(k'+k+1)^{1/2}}\right)^j
J_{k'-k+j}\bigl(2\rho\sqrt{k'+k+1}\bigr),
\end{equation}
where
\[
F_{k',k}\,:=\,\frac{k'!}{k!}\Bigl(\frac2{k'+k+1}\Bigr)^{k'-k}.
\]
The next two results give estimates for the constants appearing in
\eqref{efeimxpkpn:eq}.

\begin{lem}
\label{Ajest:lem}
Suppose $k'\ge2$ and $0\le k'-k\le k'^{\,2/3}$. Then
\[
\abs{A_j}\,\le\,(k'+k+1)^{j/3}.
\]
\end{lem}

\begin{proof}
Set $m=k'-k$ and $n=k'+k+1$ so 
\[
0\,\le\,m
\,\le\,k'^{\,2/3}
\,\le\,(k'+k+1)^{2/3}\,=\,n^{2/3}
\]
while $k'\ge2$ and $k\ge0$ so $n\ge 3$. 

\bigskip

We have $A_0=1=n^0$, $A_1=0\le n^{1/3}$ and
$m,1\le n^{2/3}$ so
$A_2=\tfrac12(m+1)\le n^{2/3}$. Now let $J\ge2$ and suppose
the result hold for $j\le J$. Since
\[
A_{J+1}\;=\;\frac{J+m}{J+1}A_{J-1}-\frac n{J+1}A_{J-2}
\]
we then get
\[
\abs{A_{J+1}}
\;\le\;\frac{J+m}{J+1}\vtsp n^{(J-1)/3}
+\frac n{J+1}\vtsp n^{(J-2)/3}
\;=\;n^{(J+1)/3}\vtsp\frac{(J+m)n^{-2/3}+1}{J+1}.
\]
Now $mn^{-2/3}\le1$ while
\begin{eqnarray*}
n\ge3
&\Longrightarrow&n^{-2/3}\le3^{-2/3}\le\tfrac12\\
&\Longrightarrow&J(1-n^{-2/3})\ge1\qquad\mbox{(as $J\ge2$)}\\
&\Longrightarrow&1+J\vtsp n^{-2/3}\le J.
\end{eqnarray*}
Thus $(J+m)\vtsp n^{-2/3}+1\le J+1$. Therefore
$\abs{A_{J+1}}\le n^{(J+1)/3}$ and the result follows by induction.
\end{proof}

\begin{lem}
\label{Fnkest:lem}
If $0\le k\le k'$ then $F_{k',k}\le 1$.
\end{lem}

\begin{proof}
We have
\[
F_{k',k}=\frac{k'(k'-1)\dots(k+1)}{\tfrac12(k'+k+1)\dots\tfrac12(k'+k+1)},
\]
where the numerator and denominator both contain $k'-k$ terms. Now set
$m=\frac12(k'-k-1)$ and $n=\frac12(k'+k+1)$ so $m\le n$ while 
\[
F_{k',k}=\frac{(n+m)}n\frac{(n+m-1)}n\dots\frac{(n-m-1)}n\frac{(n-m)}n.
\]
If $k'-k$ is odd this can be rearranged as
\[
F_{k',k}=\frac{(n+m)(n-m)}{n^2}\frac{(n+m-1)(n-m-1)}{n^2}\dots\frac nn,
\]
while if $k'-k$ is even we get
\[
F_{k',k}=\frac{(n+m)(n-m)}{n^2}\frac{(n+m-1)(n-m-1)}{n^2}\dots
\frac{(n+\tfrac12)(n-\tfrac12)}{n^2}.
\]
The result now follows from the fact that
\[
\frac{(n+m')(n-m')}{n^2}\,=\,\frac{n^2-m'^{\,2}}{n^2}\,\le\,1
\]
for any $0\le m'\le n$. 
\end{proof}

Next we obtain some estimates for the Bessel functions appearing in
\eqref{efeimxpkpn:eq}.  

\begin{lem}
\label{uglyest:lem}
For any $x,\epsilon>0$ and $n\in[0,x/2]$
\[
\bigabs{\bigl\{\theta\in[0,\pi]\,\big|\,
\abs{x\cos(\theta)-n}<\epsilon\bigr\}}
\;\le\;\frac{4\pi}3\vtsp\frac\epsilon x.
\]
\end{lem}

\begin{proof}
Set $\delta=\epsilon/x$, $y=n/x$ and 
$\Omega_{y,\delta}=\Cos^{-1}([y-\delta,y+\delta])$; we need to show
that $\abs{\Omega_{y,\delta}}\le4\pi\delta/3$. 

Now set $\theta_0=\Cos^{-1}(y)$ and let $\ell(\theta)$ denote the affine
function with $\ell(0)=1$ and $\ell(\theta_0)=y$. It is easy to see
that $\abs{\cos(\theta)-y}\ge\abs{\ell(\theta)-y}$ which implies
$\abs{\Omega_{y,\delta}}\le2\delta/\abs{L}$ where $L$ is the gradient of
$\ell(\theta)$. On the other hand, $y\in[0,\frac12]$ so the minimum
value for $\abs{L}$ occurs when $y=1/2$; hence
$1/\abs{L}\le2\Cos^{-1}(1/2)=2\pi/3$ and the result follows.
\end{proof}

\begin{lem}
\label{Besselest:lem}
For any $n\in\NZ$ and $x\ge2n$ we have $\abs{J_n(x)}\le4x^{-1/2}$.
\end{lem}

Surely this estimate (or an improvement) lies in a book somewhere!

\begin{proof}
Define a function by $f(\theta)=x\sin(\theta)-n\theta$ so we have the
following integral representation for the Bessel function $J_n$ (see
9.1.21 in \cite{AS}); 
\begin{equation}
\label{intformofJn:eq}
J_n(x)\,=\,\frac1\pi\int_0^\pi\cos(f(\theta))\d\theta.
\end{equation}
Now set
\[
\Omega_0=\bigl\{\theta\in[0,\pi]\,\big|\,
\abs{f'(\theta)}<x^{1/2}\bigr\}
\quad\text{and}\quad
\Omega_1=[0,\pi]\!\setminus\!\Omega_0
\]
so $J_n(x)=(I_0+I_1)/\pi$ where 
$\,I_k=\int_{\Omega_k}\!\cos(f(\theta))\d\theta\,$
for $k=0,1$. Lemma \ref{uglyest:lem} gives
\begin{equation}
\label{estforI0:eq}
\abs{I_0}\;\le\;\abs{\Omega_0}
\;\le\;\frac{4\pi}3x^{-1/2}.
\end{equation}
On the other hand 
\[
I_1\;=\;\biggl[\frac{\sin(f(\theta))}%
  {f'(\theta)}\biggr]_{\partial\Omega_1}
+\int_{\Omega_1}\frac{f''(\theta)}%
  {(f'(\theta))^2}\sin(f(\theta))\d\theta.
\]
Now $f''(\theta)=-x\sin(\theta)\le0$ on $[0,\pi]$ while
$(f'(\theta))^2>0$ on $\Omega_1$. Thus
\[
\lrabs{\int_{\Omega_1}\frac{f''(\theta)}%
  {(f'(\theta))^2}\sin(f(\theta))\d\theta}
\;\le\;-\int_{\Omega_1}\frac{f''(\theta)}%
  {(f'(\theta))^2}\d\theta
\;=\;\biggl[\frac{1}%
  {f'(\theta)}\biggr]_{\partial\Omega_1}.
\]
Furthermore $f'(\theta)$ is decreasing on $[0,\pi]$ so $\Omega_0$ 
consists of a single interval. Hence
$\partial\Omega_1\!\setminus\!\{0,\pi\}$ contains at most 2
points. Since $f(0)=0$ and $f(\pi)=-n\pi$ we then get
\begin{equation}
\label{estforI1:eq}
\abs{I_1}\;\le\;
\lrabs{\biggl[\frac{\sin(f(\theta))}%
  {f'(\theta)}\biggr]_{\partial\Omega_1}}
\,+\,\biggl[\frac{1}{f'(\theta)}\biggr]_{\partial\Omega_1}
\;\le\;6\max_{\theta\in\Omega_1}\frac1{\abs{f'(\theta)}}
\;\le\;6x^{-1/2}.
\end{equation}
Combining \eqref{estforI0:eq}, \eqref{estforI1:eq} we now get
\[
\abs{J_n(x)}\;\le\;\frac1\pi\bigl(\abs{I_0}+\abs{I_1}\bigr)
\;\le\;\frac1\pi\Bigl(\frac{4\pi}3+6\Bigr)x^{-1/2}
\;\le\;4x^{-1/2},
\]
completing the result.
\end{proof}

\begin{lem}
\label{matest:lem}
Suppose $k'\ge2$, $0\le k'-k\le\rho(k'+k+1)^{1/2}$ and $2\rho\le(k'+k+1)^{1/6}$. Then
\[
\abs{\Ucpt{k,k'}}\,\le\,
\bigl(4(2\rho)^{-1/2}+\tfrac12(2\rho)^2\bigr)\vtsp(k'+k+1)^{-1/4}.
\]
\end{lem}

Before starting, note that as a clear consequence of
\eqref{intformofJn:eq} we have
\begin{equation}
\label{basicbesselest:eq}
\abs{J_n(x)}\,\le\,1.
\end{equation}

\begin{proof}
Since $2,k\le k'$
\[
k'-k\;\le\;\frac12(k'+k+1)^{2/3}\;\le\;\frac12\Bigl(\frac52\Bigr)^{2/3}k'^{2/3}
\;\le\;k'^{2/3}.
\]
Now combining \eqref{efeimxpkpn:eq} with \eqref{defAja:eq},
\eqref{basicbesselest:eq} and Lemmas \ref{Ajest:lem} and
\ref{Fnkest:lem} we get
\begin{align*}
\abs{\Ucpt{k,k'}}\;&\le\;\sqrt{\smash[b]{\big.F_{k',k}}}\,
\sum_{j=0}^\infty\,\abs{A_j}\,
\frac{\rho^j}{(k'+k+1)^{j/2}}\,
\bigabs{J_{k'-k+j}\bigl(2\rho\sqrt{k'+k+1}\bigr)}\\
&\le\;\bigabs{J_{k'-k}\bigl(2\rho\sqrt{k'+k+1}\bigr)}+
\sum_{j\ge2}^\infty\rho^j(k'+k+1)^{-j/6}\\
&\le\;\bigabs{J_{k'-k}\bigl(2\rho\sqrt{k'+k+1}\bigr)}+
\tfrac12(2\rho)^2(k'+k+1)^{-1/3},
\end{align*}
where the last line follows from the hypothesis that
$\rho(k'+k+1)^{-1/6}\le1/2$. Lemma \ref{Besselest:lem} can now be used
to estimate the remaining Bessel function term.
\end{proof}

\subsection*{Main estimate}

The next result is the main estimate we will need for the matrix
elements $\abs{\ipd{\pot\efn{k}}{\efn{k'}}}$. This estimate is valid
in a parabolic region around the diagonal $k=k'$; the width of this
region is governed by the quantity
\[
\md\,:=\,\min_{\a\in\qlat}\,\aba\a\vtsp,
\]
which is positive since $\qlat$ is discrete and doesn't contain
$\mathbf 0$. Although not required in this paper, we remark that
for a general parabolic region around the diagonal one is restricted
to estimates of the form 
$\abs{\ipd{\pot\efn{k}}{\efn{k'}}}\le C(V)(k'+k+1)^{-1/6}$.

\begin{prop}
\label{keyest:prop}
Suppose $\pot$ satisfies condition \eqref{potcond1:eq} and set 
\begin{equation}
\label{defofkappa:eq}
\kappa\,=\,\min\{1/3,\md/(2\sqrt3)\}.
\end{equation}
If $n\in\N$ and $k,k'\in\NZ$ satisfy
$\abs{k-n},\abs{k'-n}\le\kappa n^{1/2}$ then
\begin{equation}
\label{keyest:eq}
\abs{\ipd{\pot\efn{k}}{\efn{k'}}}
\,\le\,C(\pot)\vtsp n^{-1/4}.
\end{equation}
\end{prop}

\begin{proof}
We have $\abs{\ipd{\pot\efn{k}}{\efn{k'}}}\le\norm{V}$ for any
$k,k'\in\NZ$ so we can increase $C(V)$ if necessary to ensure that 
\eqref{keyest:eq} is satisfied for $n=1,2$. Furthermore $\pot$ is
self-adjoint so
$\abs{\ipd{\pot\efn{k'}}{\efn{k}}}=\abs{\ipd{\pot\efn{k}}{\efn{k'}}}$. 
It thus suffices to prove the result assuming $n\ge3$ and $k',k\in\NZ$
satisfy $k'\ge k$ and $\abs{k-n},\abs{k'-n}\le\kappa n^{1/2}$. Then
$k',k\ge n-\frac13n^{1/2}\ge\frac23n$ so $k'\ge2$, 
\begin{equation}
\label{ineqtorelback1:eq}
k'+k+1\,\ge\,\frac43n
\end{equation}
and
\begin{equation}
\label{ineqtorelback2:eq}
0\,\le\,k'-k\,\le\,2\kappa n^{1/2}
\,\le\,\frac\md2(k'+k+1)^{1/2}.
\end{equation}
Now set $K=(k'+k+1)^{1/6}$. Using \eqref{defofpot:eq},
\eqref{defofUcpt:eq} and \eqref{UisUso1:eq} we have
\begin{equation}
\label{estforpotl1:eq}
\abs{\ipd{\pot\efn{k}}{\efn{k'}}}
\;\le\;\sum_{\a\in\qlat}\,\abs{\Ucpt{k,k'}}\vtsp\abs{\potcpt}
\;\le\;\sum_{\substack{\a\in\qlat\\\aba\a\le K}}
\abs{\Ucpt{k,k'}}\vtsp\abs{\potcpt}
+\sum_{\substack{\a\in\qlat\\\aba\a>K}}\abs{\potcpt}.
\end{equation}
Since $1<K^{-3/2}\aba\a^{3/2}$ whenever $\aba\a>K$,
\eqref{potcond1:eq} and \eqref{ineqtorelback1:eq} give us
\[
\sum_{\substack{\a\in\qlat\\\aba\a>K}}\abs{\potcpt}
\;\le\;K^{-3/2}\sum_{\a\in\qlat}\,\aba\a^{3/2}\abs{\potcpt}
\;\le\;C(\pot)\vtsp n^{-1/4}.
\]
Now let $\a\in\qlat$. Since $\aba\a=2\rho$ (see \eqref{defofrho:eq})
the definition of $\md$ implies $\md/2\le\rho$ and thus
$k'-k\le\rho K^3$ by \eqref{ineqtorelback2:eq}. 
Lemma \ref{matest:lem}, \eqref{potcond1:eq} and
\eqref{ineqtorelback1:eq} then give
\[
\sum_{\substack{\a\in\qlat\\\aba\a\le K}}
\abs{\Ucpt{k,k'}}\vtsp\abs{\potcpt}
\;\le\;K^{-3/2}\sum_{\a\in\qlat}\,(4\aba\a^{-1/2}+\tfrac12\aba\a^2)
\vtsp\abs{\potcpt}
\;\le\;C(\pot)\vtsp n^{-1/4}.
\]
The result follows.
\end{proof}

\subsection*{First order term}

The next result is used to obtain the explicit form for the first
order correction term in the asymptotics for $\lambda_n(\op+\pot)$. 

\begin{prop}
\label{1storderterm:prop}
Suppose $\pot$ satisfies condition \eqref{potcond1:eq}. 
Then 
\[
\ipd{\pot\efn{n}}{\efn{n}}
\,=\,W(\sqrt{n})\vtsp n^{-1/4}+O(n^{-1/2})
\]
as $n\to+\infty$, where $W$ is defined by \eqref{defofauxfnU:eq}.
\end{prop}

\begin{proof}
Let $\a\in\TR\!\setminus\!\{\mathbf 0\}$ and set $\rho=\aba\a/2$. 
Using \eqref{efeimxpkpn:eq} and the fact that $F_{n,n}=1$ we get
\[
\Ucpt{n,n}\,=\,\sum_{j=0}^\infty\,A_j\vtsp\frac{\rho^j}{(2n+1)^{j/2}}
\vtsp J_j\bigl(2\rho\sqrt{2n+1}\bigr).
\]
Now suppose $2\rho\le N$ where $N:=(2n+1)^{1/6}$. Using
\eqref{defAja:eq}, \eqref{basicbesselest:eq} and Lemma \ref{Ajest:lem} we have
\begin{align*}
\bigabs{\Ucpt{n,n}-J_0(2\rho\sqrt{2n+1})}
\;&\le\;\tfrac12\rho^2(2n+1)^{-1}+
\sum_{j\ge3}^\infty\rho^j(2n+1)^{-j/6}\\
&\le\;\tfrac18\aba\a^2(2n+1)^{-1}+\tfrac14\aba\a^3(2n+1)^{-1/2}.
\end{align*}
Standard asymptotic forms for Bessel functions (see 9.2.1 in
\cite{AS}) give us
\[
J_0(z)
\;=\;\sqrt{\frac{2}{\pi z}}\cos\Bigl(z-\frac\pi4\Bigr)+O(z^{-3/2})
\]
while 
\[
\lrabs{\frac{d}{dz}\Bigl(\frac1{\sqrt z}\cos\Bigl(z-\frac\pi4\Bigr)\Bigr)}
\,\le\,z^{-1/2}+\frac12z^{-3/2}
\]
and $2\rho\sqrt{2n+1}-2\rho\sqrt{2n}\le2^{-1/2}\rho n^{-1/2}$. 
It follows that
\begin{eqnarray*}
\lefteqn{
\lrabs{J_0(2\rho\sqrt{2n+1})-\sqrt{\frac2\pi}(2\rho)^{-1/2}
(2n)^{-1/4}\cos\Bigl(2\rho\sqrt{2n}-\frac\pi4\Bigr)}}\qquad\\
&\le&C(2\rho)^{-3/2}(2n+1)^{-3/4}\\
&&\qquad{}+\sqrt{\frac2\pi}\Bigl((2\rho)^{-1/2}(2n)^{-1/4}
+\frac12(2\rho)^{-3/2}(2n)^{-3/4}\Bigr)
\vtsp2^{-1/2}\rho\vtsp n^{-1/2}\\
&\le&C((2\rho)^{-3/2}+(2\rho)^{1/2})n^{-3/4}.
\end{eqnarray*}
Combining the above estimates we thus obtain
\[
\lrabs{\Ucpt{n,n}-\frac{2^{1/4}}{\sqrt\pi}\aba\a^{-1/2}n^{-1/4}
\cos\Bigl(\aba\a\sqrt{2n}-\frac\pi4\Bigr)}
\;\le\;C(\aba\a^{-3/2}+\aba\a^3)n^{-1/2}
\]
whenever $\aba\a\le N$. Using \eqref{defofpot:eq}, \eqref{defofauxfnU:eq},
\eqref{defofUcpt:eq} and \eqref{UisUso1:eq} we thus have
\begin{eqnarray*}
\lefteqn{
\bigabs{\ipd{\pot\efn{n}}{\efn{n}}-W(\sqrt{n})\vtsp n^{-1/4}}}\qquad\\
&\le&Cn^{-1/2}\sum_{\substack{\a\in\qlat\\\aba\a\le N}}
(\aba\a^{-3/2}+\aba\a^{3})\vtsp\abs{\potcpt}
+\sum_{\substack{\a\in\qlat\\\aba\a>N}}(1+\aba\a^{-1/2})\abs{\potcpt}.
\end{eqnarray*}
Since $1<N^{-3}\aba\a^3<n^{-1/2}\aba\a^3$ whenever $\aba\a>N$ the term inside
the last sum can be replaced with
$n^{-1/2}(\aba\a^3+\aba\a^{5/2})\abs\potcpt$. Using
\eqref{potcond1:eq} we then get
\[
\bigabs{\ipd{\pot\efn{n}}{\efn{n}}-W(\sqrt{n})\vtsp n^{-1/4}}
\;\le\;Cn^{-1/2}\sum_{\a\in\qlat}\,(\aba\a^{-3/2}+\aba\a^3)\vtsp\abs{\potcpt}
\;\le\;C(V)n^{-1/2},
\]
completing the result.
\end{proof}

\section{Resolvent Estimates}
\label{resest:sec}
\newcommand{\cont}[1][n]{\Gamma_{\epsilon,{#1}}}

For any $\lambda\in\C\!\setminus\!\sigma(\op)$ let
$\res(\lambda)\,=\,(\op-\lambda)^{-1}$
denote the resolvent of the operator $\op$; we will also write $\res$
for $\res(\lambda)$ where this should not cause confusion.

Let $\kappa$ denote the constant defined in \eqref{defofkappa:eq}.
For a given $n\in\N$ we will make repeated use of the 
partition of $\NZ$ defined by 
\begin{equation}
\label{genpart:eq}
I=\bigl\{k\in\NZ\,\big|\,\abs{k-n}\le\kappa n^{1/2}\bigr\}
\quad\mbox{and}\quad
J=\NZ\!\setminus\!I.
\end{equation}
For any
$\epsilon\in(0,\alpha)$ and $n\in\NZ$, let $\cont$ be the anti-clockwise
circular contour in $\C$ centred at
$\lambda_n=\lambda_n(\op)=\alpha(2n+1)$. If $\lambda\in\cont$ then
$\lambda=\alpha(2n\!+\!1)+\epsilon e^{i\theta}$ for some
$\theta\in[0,2\pi)$. It follows that  
$\abs{\lambda-\lambda_k}=\abs{2\alpha(n\!-\!k)+\epsilon e^{i\theta}}$
for any $k\in\NZ$. Straightforward arguments then lead to the
following estimates;
\begin{equation}
\label{epsestlam1:eq}
\sum_{k\in I}\,\abs{\lambda-\lambda_k}^{-1}
\,\le\,C(\epsilon)\,\ln(n),
\end{equation}
\begin{equation}
\label{epsestlam2a:eq}
\sum_{k\in\NZ}\,\abs{\lambda-\lambda_k}^{-2}
\,\le\,C(\epsilon),
\end{equation}
\begin{equation}
\label{epsestlam2:eq}
\sum_{k\in J}\,\abs{\lambda-\lambda_k}^{-2}
\,\le\,C\vtsp n^{-1/2}
\end{equation}
and
\begin{equation}
\label{epsestlam3:eq}
\abs{\lambda-\lambda_k}\ge C\vtsp n^{1/2}
\quad\text{for any $k\in J$.}
\end{equation}

\bigskip

The first two results in this section relate to the operator
$\res(\lambda)\pot\res(\lambda)$, which is 
clearly bounded whenever $\lambda$ is in the resolvent set of $\op$. 
We show that it is in fact trace class while its
operator norm decreases as $n^{-1/4}$ for $\lambda\in\cont$.

\begin{lem}
\label{decbndonR0VR0:lem}
For any $n\in\N$ and $\lambda\in\cont$ we have
\[
\norm{\res(\lambda)\pot\res(\lambda)}
\,\le\,\norm{\res(\lambda)\pot\res(\lambda)}_2
\,\le\,C(\pot,\epsilon)\vtsp n^{-1/4}.
\]
\end{lem}

We remark that since
$\{\phi_k\,|\,k\in\NZ\}$ is an orthonormal basis of $L^2(\R)$ 
\begin{equation}
\label{phicompVbnd:eq}
\sum_{k'\in\NZ}\,\abs{\ipd{\pot\phi_k}{\phi_{k'}}}^2
\;=\;\norm{\pot\phi_k}^2
\;\le\;\norm{\pot}^2.
\end{equation}

\begin{proof}
Using the orthonormal basis $\{\efn{k}\,|\,k\in\NZ\}$ we have
\begin{equation}
\label{HSnormofRVR:eq}
\norm{\res\pot\res}^2_2
\;=\;\sum_{k,k'\in\NZ}
\abs{\ipd{\res\pot\res\efn{k}}{\efn{k'}}}^2\\
\;=\;\sum_{k,k'\in\NZ}
\frac{\abs{\ipd{\pot\efn{k}}{\efn{k'}}}^2}%
 {\abs{\lambda_k-\lambda}^2\abs{\lambda_{k'}-\lambda}^2}.
\end{equation}
We will split this sum using the partition \eqref{genpart:eq}. Firstly
Proposition \ref{keyest:prop} and \eqref{epsestlam2a:eq} imply
\[
\sum_{k,k'\in I}\frac{\abs{\ipd{\pot\efn{k}}{\efn{k'}}}^2}%
 {\abs{\lambda_k-\lambda}^2\abs{\lambda_{k'}-\lambda}^2}
\;\le\;C(\pot)n^{-1/2}\biggl(\,\sum_{k\in I}
\frac1{\abs{\lambda_k-\lambda}^2}\biggr)\Big.^2
\;\le\;C(\pot,\epsilon)n^{-1/2}.
\]
Now using \eqref{epsestlam3:eq}, \eqref{phicompVbnd:eq} and
\eqref{epsestlam2a:eq} we get
\begin{align*}
\sum_{\substack{k\in\NZ\\k'\in J}}
\frac{\abs{\ipd{\pot\efn{k}}{\efn{k'}}}^2}%
 {\abs{\lambda_k-\lambda}^2\abs{\lambda_{k'}-\lambda}^2}
\;&\le\;Cn^{-1}\sum_{k\in\NZ}\frac1{\abs{\lambda_k-\lambda}^2}
\sum_{k'\in J}\,\abs{\ipd{\pot\efn{k}}{\efn{k'}}}^2\\
&\le\;C(\pot,\epsilon)n^{-1}.
\end{align*}
The remaining part of the sum on the right hand side of
\eqref{HSnormofRVR:eq} involves $k\in J$ and $k'\in I\subset\NZ$; thus
we can estimate this part using an argument similar to the last one with
$k$ and $k'$ swapped. 
\end{proof}

\begin{lem}
\label{traceestR0VR0:lem}
For any $n\in\NZ$ and $\lambda\in\cont$ the operator
$\res(\lambda)\pot\res(\lambda)$ is trace class. Furthermore
$\norm{\res(\lambda)\pot\res(\lambda)}_1$ is uniformly bounded (in
$n$ and $\lambda\in\cont$).
\end{lem}

\begin{proof}
The set $\{\efn{k}\,|\,k\in\NZ\}$ is an orthonormal eigenbasis
for $\res$ with corresponding eigenvalues
$(\lambda_k-\lambda)^{-1}$, $k\in\NZ$ so \eqref{epsestlam2a:eq}
implies 
\[
\norm{\res}_2^2
\,=\,\sum_{k\in\NZ}\,\abs{\lambda_k-\lambda}^{-2}
\,\le\,C(\epsilon).
\]
Thus
$\norm{\res\pot\res}_1\,=\,\norm{\pot\res^2}_1
\,\le\,\norm{\pot}\,\norm{\res^2}_1
\,\le\,\norm{\pot}\,\norm{\res}_2^2
\,\le\,C(\epsilon)\vtsp\norm{\pot}$.
\end{proof}

\medskip

Suppose $n\in\NZ$ and $j\in\N$. From the previous result we know that
$\res(\lambda)\pot\res(\lambda)$ is trace class for any
$\lambda\in\cont$. On the other hand $\res(\lambda)\pot$ is bounded
(in fact $\norm{\res(\lambda)\pot}\le\epsilon^{-1}\norm{\pot}$). It follows
that
\[
(\res(\lambda)\pot)^j\res(\lambda)
\,=\,(\res(\lambda)\pot)^{j-1}\res(\lambda)\pot\res(\lambda)
\]
is also trace class with trace norm uniformly bounded for
$\lambda\in\cont$. The work in the remainder of this section leads to
Proposition \ref{estonTrk:prop} where we obtain an estimate for the
trace of an integral of such operators.

\begin{lem}
\label{ftogindst:lem}
Let $n\ge2$, $\lambda\in\cont$ and suppose $f:\NZ\to\C$ satisfies
\begin{equation}
\label{indhypa:eq}
\sum_{k\in\NZ}\,\abs{f(k)}^2
\,\le\,C_1^2
\quad\text{and}\quad
\abs{f(k)}\,\le\,C_1n^{-1/4}
\ \text{when $k\in I$}
\end{equation}
for some constant $C_1$. For each $k\in\NZ$ set
\[
g(k)=\sum_{k'\in\NZ}\,\frac{f(k')\,\ipd{\pot\efn{k'}}{\efn{k}}}{\lambda-\lambda_{k'}}.
\]
Then there exists a constant $K=K(\pot,\epsilon)$ such that
\[
\sum_{k\in\NZ}\,\abs{g(k)}^2
\,\le\,C_1^2K^2n^{-1/2}\ln^2(n)
\quad\text{and}\quad
\abs{g(k)}\,\le\,C_1Kn^{-1/2}\ln(n)
\ \text{when $k\in I$.}
\]
\end{lem}

\begin{proof}
Since
\[
\sum_{k\in\NZ}\,\abs{g(k)}^2
\;=\;\sum_{k,k',k''\in\NZ}\frac{f(k')\overline{f(k'')}%
  \ipd{\pot\efn{k'}}{\efn{k}}\ipd{\efn{k}}{\pot\efn{k''}}}%
  {(\lambda-\lambda_{k'})\overline{(\lambda-\lambda_{k''})}}
\]
and
\[
\biggl\lvert\sum_{k\in\NZ}
\,\ipd{\pot\efn{k'}}{\efn{k}}\ipd{\efn{k}}{\pot\efn{k''}}\biggr\rvert
\;=\;\bigabs{\ipd{\pot\efn{k'}}{\pot\efn{k''}}}
\;\le\;\norm{\pot}^2
\]
it follows that 
\[
\sum_{k\in\NZ}\,\abs{g(k)}^2
\;\le\;\norm{\pot}^2\biggl(\sum_{k\in I}\,
\frac{\abs{f(k)}}{\abs{\lambda-\lambda_k}}
\,+\,\sum_{k\in J}\,\frac{\abs{f(k)}}{\abs{\lambda-\lambda_k}}\biggr)\Big.^2.
\]
Using the second part of \eqref{indhypa:eq} and \eqref{epsestlam1:eq} we get
\[
\sum_{k\in I}\,
\frac{\abs{f(k)}}{\abs{\lambda-\lambda_k}}
\;\le\;C_1n^{-1/4}\sum_{k\in I}\,\abs{\lambda-\lambda_k}^{-1}
\;\le\;C_1C_2(\epsilon)n^{-1/4}\ln(n).
\]
On the other hand the first part of \eqref{indhypa:eq} and \eqref{epsestlam2:eq} give
\begin{eqnarray*}
\sum_{k\in J}\,
\frac{\abs{f(k)}}{\abs{\lambda-\lambda_k}}
&\le&\biggl(\sum_{k\in J}\,\abs{f(k)}^2\biggr)\Big.^{1/2}
\biggl(\sum_{k\in J}\,\abs{\lambda-\lambda_k}^{-2}\biggr)\Big.^{1/2}\\
&\le&C_1C_3n^{-1/4}
\ \le\ 2C_1C_3n^{-1/4}\ln(n).
\end{eqnarray*}
Putting these estimates together now leads to 
\[
\sum_{k\in\NZ}\,\abs{g(k)}^2\,\le\,C_1^2K_1^2 n^{-1/2}\ln^2(n),
\]
with $K_1=\norm{\pot}(C_2(\epsilon)+2C_3)$.
Now suppose $k\in I$ and write $g(k)=g_I(k)+g_J(k)$ where
\[
g_I(k)=\sum_{k'\in I}\,\frac{f(k')\,\ipd{\pot\efn{k'}}{\efn{k}}}%
  {\lambda-\lambda_{k'}}
\quad\mbox{and}\quad
g_J(k)=\sum_{k'\in J}\,\frac{f(k')\,\ipd{\pot\efn{k'}}{\efn{k}}}%
  {\lambda-\lambda_{k'}}.
\]
From Proposition \ref{keyest:prop}, the second part of
\eqref{indhypa:eq} and \eqref{epsestlam1:eq} we get
\[
\abs{g_I(k)}
\;\le\;C_1C(\pot)n^{-1/2}\sum_{k'\in I}\,\abs{\lambda-\lambda_{k'}}^{-1}
\;\le\;C_1C_4(\pot,\epsilon)\vtsp n^{-1/2}\ln(n).
\]
On the other hand \eqref{epsestlam3:eq}, the first part of
\eqref{indhypa:eq} and \eqref{phicompVbnd:eq} give us
\begin{eqnarray*}
\abs{g_J(k)}&\le&C(\epsilon)n^{-1/2}
\biggl(\sum_{k'\in J}\,\abs{f(k')}^2\biggr)\Big.^{1/2}
\biggl(\sum_{k'\in J}\,\bigabs{\ipd{\pot\efn{k'}}{\efn{k}}}^2\biggr)\Big.^{1/2}\\
&\le&C_1C_5(\epsilon)\norm{\pot}n^{-1/2}
\ \le\ 2C_1C_5(\epsilon)\norm{\pot}n^{-1/2}\ln(n).
\end{eqnarray*}
Putting these estimates together now leads to 
$\abs{g(k)}\le C_1K_2n^{-1/2}\ln(n)$ with 
$K_2=C_4(\pot,\epsilon)+2C_5(\epsilon)\norm{\pot}$.
Taking $K=\max\{K_1,K_2\}$, completes the result.
\end{proof}

Taking $f(k)=\ipd{\pot\efn{n}}{\efn{k}}$ we can use
\eqref{phicompVbnd:eq} and Proposition \ref{keyest:prop} to check that
\eqref{indhypa:eq} is satisfied. The next result then follows from Lemma
\ref{ftogindst:lem} by use of induction; we can take
$K=\max\{\norm\pot,C(\pot),K'\}$ where $C(\pot)$ and $K'$ are the
constants coming from Proposition \ref{keyest:prop} and Lemma
\ref{ftogindst:lem} respectively. 

\begin{lem}
\label{bndonintgrand:lem}
Suppose $n\ge2$ and $j\in\NZ$. Then there exists a constant
$K=K(\pot,\epsilon)$ such that for all $\lambda\in\cont$ we have
\[
\Biggl\lvert\sum_{k_1,\dots,k_j\in\NZ}
\frac{\ipd{\pot\efn{n}}{\efn{k_1}}\ipd{\pot\efn{k_1}}{\efn{k_2}}
  \dots\ipd{\pot\efn{k_j}}{\efn{n}}}%
  {(\lambda-\lambda_{k_1})\dots(\lambda-\lambda_{k_j})}\Biggr\rvert
\;\le\;K^{j+1}n^{-(j+1)/4}\ln^j(n).
\]
\end{lem}

\begin{prop}
\label{estonTrk:prop}
Suppose $n\ge2$ and $j\in\N$. Then
\[
\lrabs{\Trace\frac1{2\pi i}\oint_{\cont}
\lambda\,\res(\lambda)(\pot\res(\lambda))^j\d\lambda}
\;\le\;K^jn^{-j/4}\ln^{j-1}(n),
\]
where $K$ is the constant from Lemma \ref{bndonintgrand:lem}.
\end{prop}

\begin{proof}
Since $\{\efn{k'}\,|\,k'\in\NZ\}$ is an orthonormal basis of $L^2(\R)$ we
have
\[
(\pot\res)\efn{k}
\;=\;\sum_{k'\in\NZ}\ipd{\pot\res\efn{k}}{\efn{k'}}\efn{k'}
\;=\;\sum_{k'\in\NZ}\frac{\ipd{\pot\efn{k}}{\efn{k'}}}{\lambda_k-\lambda}\efn{k'}.
\]
Continuing by induction we get
\[
(\pot\res)^j\efn{k}
\;=\;\sum_{k_1,\dots,k_j\in\NZ}
\frac{\ipd{\pot\efn{k}}{\efn{k_1}}\ipd{\pot\efn{k_1}}{\efn{k_2}}
  \dots\ipd{\pot\efn{k_{j-1}}}{\efn{k_j}}}
  {(\lambda_k-\lambda)(\lambda_{k_1}-\lambda)\dots(\lambda_{k_{j-1}}-\lambda)}
\efn{k_j}.
\]
Together with the fact that 
$\ipd{\res\efn{k_j}}{\efn{k_0}}=\delta_{k_j,k_0}(\lambda_{k_0}-\lambda)^{-1}$
we now get
\begin{eqnarray*}
\lefteqn{
\Trace(\res(\pot\res)^j)
\ =\ \sum_{k_0\in\NZ}\,\bigipd{\res(\pot\res)^j\phi_{k_0}}{\phi_{k_0}}}\quad\\
&=&\sum_{k_0,\dots,k_j\in\NZ}
\frac{\ipd{\pot\efn{k_0}}{\efn{k_1}}\ipd{\pot\efn{k_1}}{\efn{k_2}}
  \dots\ipd{\pot\efn{k_{j-1}}}{\efn{k_j}}}
  {(\lambda_{k_0}-\lambda)(\lambda_{k_1}-\lambda)\dots(\lambda_{k_{j-1}}-\lambda)}
\ipd{\res\efn{k_j}}{\efn{k_0}}\\
&=&\sum_{k_0,\dots,k_{j-1}\in\NZ}
\frac{\ipd{\pot\efn{k_0}}{\efn{k_1}}\ipd{\pot\efn{k_1}}{\efn{k_2}}
  \dots\ipd{\pot\efn{k_{j-1}}}{\efn{k_0}}}
  {(\lambda_{k_0}-\lambda)^2(\lambda_{k_1}-\lambda)\dots(\lambda_{k_{j-1}}-\lambda)}
\ =\ \sum_{l=0}^{j-1}\,\frac1{\lambda_{k_l}-\lambda}A(\lambda),
\end{eqnarray*}
where $A(\lambda)$ is the meromorphic function
\begin{equation}
\label{defofA:eq}
A(\lambda)\;=\;\frac1j\sum_{k_0,\dots,k_{j-1}\in\NZ}
\frac{\ipd{\pot\efn{k_0}}{\efn{k_1}}\ipd{\pot\efn{k_1}}{\efn{k_2}}
  \dots\ipd{\pot\efn{k_{j-1}}}{\efn{k_0}}}
  {(\lambda_{k_0}-\lambda)(\lambda_{k_1}-\lambda)\dots(\lambda_{k_{j-1}}-\lambda)}.
\end{equation}
Since
\[
\frac{d}{d\lambda}\,\lambda
  \biggl(\,\prod_{i=0}^{j-1}\,\frac1{\lambda_{k_i}-\lambda}\biggr)
\;=\;\prod_{i=0}^{j-1}\,\frac1{\lambda_{k_i}-\lambda}
\,+\,\lambda\sum_{l=0}^{j-1}\,\frac1{\lambda_{k_l}-\lambda}
  \biggl(\,\prod_{i=0}^{j-1}\,\frac1{\lambda_{k_i}-\lambda}\biggr)
\]
for any $k_0,\dots,k_{j-1}\in\NZ$, 
we can rewrite the above equation as
\[
\Trace\, \lambda \res(\pot\res)^j
\;=\;\frac{d}{d\lambda}\bigl(\lambda A(\lambda)\bigr)-A(\lambda).
\]
Integrating around the contour $\cont$ it follows that
\begin{equation}
\label{traceitoA:eq}
\Trace\,\frac1{2\pi i}\oint_{\cont}\lambda \res(\pot\res)^j\d\lambda
\;=\;-\frac1{2\pi i}\oint_{\cont}A(\lambda)\d\lambda.
\end{equation}
The poles of the meromorphic function $A(\lambda)$ occur at the
points $\lambda=\lambda_k$ for $k\in\NZ$. Since the only such point
enclosed by the contour $\cont$ is $\lambda=\lambda_n$, it follows
that the only terms in the series \eqref{defofA:eq} which contribute to the
right hand side of \eqref{traceitoA:eq} are those with at least one of
$k_0,\dots,k_{j-1}$ equal to $n$. With the help of symmetry we then
obtain the identity
\begin{eqnarray}
\lefteqn{\Trace\,\frac1{2\pi i}\oint_{\cont}\lambda
  \res(\pot\res)^j\d\lambda
}\ \nonumber\\
\label{expanTrproj:eq}
&=&-\frac1{2\pi i}\oint_{\cont}\frac1{\lambda_n-\lambda}
\sum_{k_1,\dots,k_{j-1}\in\NZ}
\frac{\ipd{\pot\phi_n}{\efn{k_1}}\ipd{\pot\efn{k_1}}{\efn{k_2}}
  \dots\ipd{\pot\efn{k_{j-1}}}{\efn{n}}}
  {(\lambda_{k_1}-\lambda)\dots(\lambda_{k_{j-1}}-\lambda)}
\d\lambda.\qquad
\end{eqnarray}
For any $\lambda\in\cont$ we have $\abs{\lambda_n-\lambda}=\epsilon$ 
while
\[
\Biggl\lvert\sum_{k_1,\dots,k_{j-1}\in\NZ}
\frac{\ipd{\pot\phi_n}{\efn{k_1}}\ipd{\pot\efn{k_1}}{\efn{k_2}}
  \dots\ipd{\pot\efn{k_{j-1}}}{\efn{n}}}
  {(\lambda_{k_1}-\lambda)\dots(\lambda_{k_{j-1}}-\lambda)}
\Biggr\rvert
\;\le\;K^jn^{-j/4}\ln^{j-1}(n)
\]
by Lemma \ref{bndonintgrand:lem}. Since the length of $\cont$ is
$2\pi\epsilon$ we finally get
\[
\lrabs{\Trace\,\frac1{2\pi i}\oint_{\cont}\lambda
  \res(\pot\res)^j\d\lambda}
\;\le\;\frac1{2\pi}\oint_{\cont}\frac1{\epsilon}\,
  K^jn^{-j/4}\ln^{j-1}(n)\d\lambda
\;=\;K^jn^{-j/4}\ln^{j-1}(n),
\]
completing the result.
\end{proof}

Taking $j=1$ in \eqref{expanTrproj:eq} leads to the formula
\begin{eqnarray}
\lefteqn{
\Trace\,\frac1{2\pi i}\oint_{\cont}\lambda
  \res(\lambda)\pot\res(\lambda)\d\lambda}
\qquad\nonumber\\
\label{formof1stTr:eq}
&&=\;-\frac1{2\pi i}\oint_{\cont}\frac1{\lambda_n-\lambda}
\ipd{\pot\efn{n}}{\efn{n}}
\d\lambda
\;=\;\ipd{\pot\efn{n}}{\efn{n}}.
\end{eqnarray}
This is needed to obtain the first order correction term in Theorem
\ref{mainres1:thm}.

\section{Proof of Theorem \ref{mainres1:thm}}
\label{mainrespf:sec}

Lemmas \ref{decbndonR0VR0:lem} and \ref{traceestR0VR0:lem} give us
$\norm{\res(\lambda)\pot\res(\lambda)}\le C_1n^{-1/4}$ and 
$\norm{\res(\lambda)\pot\res(\lambda)}_1\le C_2$ 
for all $n\in\N$ and $\lambda\in\cont$. In particular
$\norm{(\pot\res(\lambda))^2}\le\norm{\pot}C_1n^{-1/4}$.
We also note that $\norm{\res(\lambda)}=\epsilon^{-1}$.
It follows that for any $j\in\NZ$ we get
\[
\norm{(\pot\res(\lambda))^{2j}}
\;\le\;\norm{(\pot\res(\lambda))^2}^j
\;\le\;(\norm{\pot}C_1n^{-1/4})^j
\]
and
\[
\norm{(\pot\res(\lambda))^{2j+1}}
\;\le\;\norm{\pot}\norm{\res(\lambda)}\norm{(\pot\res(\lambda))^{2j}}
\;\le\;\norm{\pot}\epsilon^{-1}(\norm{\pot}C_1n^{-1/4})^j.
\]
Choose $N'\in\N$ so that $\norm{\pot}C_1{N'}^{-1/4}\le1/2$.
It follows that for any $n\ge N'$ and $\lambda\in\cont$ the series
\begin{equation}
\label{neumanseries:eq}
(I+\pot\res(\lambda))^{-1}
\,=\,\sum_{j=0}^\infty\,(-\pot\res(\lambda))^j
\end{equation}
is absolutely convergent and has norm bounded by
$2(1+\norm{\pot}\epsilon^{-1})$.
In particular, $I+\pot\res(\lambda)$ is invertible with a uniformly
bounded inverse for all $n\ge N'$ and $\lambda\in\cont$.
On the other hand, the series
\[
T(\lambda)\;:=\;\sum_{j=1}^\infty \res(\lambda)(-\pot\res(\lambda))^j
\;=\;-\res(\lambda)\pot\res(\lambda)\sum_{j=0}^\infty\,(-\pot\res(\lambda))^j
\]
is convergent in trace class with  
\[
\norm{T(\lambda)}_1
\;\le\;\norm{\res(\lambda)\pot\res(\lambda)}_1
\,\norm{(I+\pot\res(\lambda))^{-1}}
\;\le\;2C_2(1+\norm{\pot}\epsilon^{-1})
\]
for $n\ge N'$ and $\lambda\in\cont$. Setting 
\[
T_n\,=\,-\frac1{2\pi i}\oint_{\cont}\lambda T(\lambda)\d\lambda
\]
it follows that we have an absolutely convergent expansion
\[
\Trace T_n
\;=\;-\sum_{j=1}^\infty\,\Trace\frac1{2\pi i}\oint_{\cont}\lambda 
\res(\lambda)(-\pot\res(\lambda))^j\d\lambda
\]
whenever $n\ge N'$.

Now choose $N\ge N'$ so that $K N^{-1/4}\ln(N)\le1/2$
where $K$ is the constant given by Proposition \ref{estonTrk:prop}. Using this
Proposition and the above results it follows that 
\[
\lrabs{\sum_{j=2}^\infty\,\Trace\frac1{2\pi i}\oint_{\cont}\lambda 
\res(\lambda)(-\pot\res(\lambda))^j\d\lambda}
\;\le\;2K^2n^{-1/2}\ln(n)
\]
for all $n\ge N$. Therefore
\begin{align*}
\Trace T_n
\;&=\;\Trace\frac1{2\pi i}\oint_{\cont}\lambda 
\res(\lambda)\pot\res(\lambda)\d\lambda
\,+\,O(n^{-1/2}\ln(n))\\
&=\;\ipd{\pot\efn{n}}{\efn{n}}\,+\,O(n^{-1/2}\ln(n))
\end{align*}
for all $n\ge N$, where we have used \eqref{formof1stTr:eq}.

\medskip

The argument can be tied together using a standard resolvent
expansion. Set $\res_\pot(\lambda)=(\op+\pot-\lambda)^{-1}$ and let
$n\ge N$. Then
\[
\res_\pot(\lambda)
\;=\;\res(\lambda)(1+\pot\res(\lambda))^{-1}
\;=\;\res(\lambda)\sum_{j=0}^\infty\,(-\pot\res(\lambda))^j.
\]
The right hand side of
\eqref{neumanseries:eq} will still converge if $\pot$ is replaced with
$g\pot$ for some $g\in[0,1]$. Hence
$\sigma(\op+g\pot)\cap\cont=\emptyset$. Since the eigenvalues of
$\op+g\pot$ depend continuously on $g$, it follows that $\cont$ must enclose
$\lambda_n(\op+\pot)$ but no other points of
$\sigma(\op+\pot)$. Thus we can write 
%\[
%\lambda_n(\op+\pot)\,=\,-\frac1{2\pi i}\Trace\oint_{\cont}
%\lambda R_\pot(\lambda)\d\lambda
%\quad\mbox{and}\quad
%\lambda_n(\op)\,=\,-\frac1{2\pi i}\Trace\oint_{\cont}
%\lambda \res(\lambda)\d\lambda
%\]
%and so
\begin{eqnarray*}
\lefteqn{\lambda_n(\op+\pot)-\lambda_n(\op)
\;=\;-\frac1{2\pi i}\Trace\oint_{\cont}
\lambda(\res_\pot(\lambda)-\res(\lambda))\d\lambda}\qquad\qquad\\
&=&-\frac1{2\pi i}\Trace\oint_{\cont}
\lambda\sum_{j=1}^\infty 
\res(\lambda)(-\pot\res(\lambda))^j\d\lambda\\
&=&\Trace T_n
\ =\ \ipd{\pot\efn{n}}{\efn{n}}+O(n^{-1/2}\ln(n)).
\end{eqnarray*}
Theorem \ref{mainres1:thm} now follows from Proposition \ref{1storderterm:prop}.

\subsection*{Acknowledgements}

The author wishes to thank A.\ B.\ Pushnitski for several useful
discussions, especially regarding some of the special function results
used in Section \ref{matres:sec}.

\end{document}